\documentclass[12pt]{amsart}

\usepackage{amssymb}
\usepackage{amsmath}
\usepackage{amsthm}
\usepackage{graphics}

\allowdisplaybreaks

\numberwithin{equation}{section}

\theoremstyle{plain}
\newtheorem{thm}{Theorem}[section]

\newtheorem*{referthmB}{Theorem B}

\newtheorem*{referpropA}{Proposition A}

\theoremstyle{definition}

\newcommand{\ichi}{\mathbf{1}}
\newcommand{\R}{\mathbb{R}}

\newcommand{\N}{\mathbb{N}}

\newcommand{\calF}{\mathcal{F}}

\newcommand{\calS}{\mathcal{S}}
\newcommand{\supp}{\mathrm{supp}\, }

\begin{document}

\title{Bilinear pseudo-differential operators with exotic symbols, II}

\author{Akihiko Miyachi \and Naohito Tomita}
\date{}

\address{Akihiko Miyachi \\
Department of Mathematics \\
Tokyo Woman's Christian University \\
Zempukuji, Suginami-ku, Tokyo 167-8585, Japan}
\email{miyachi@lab.twcu.ac.jp}

\address{Naohito Tomita \\
Department of Mathematics \\
Graduate School of Science \\
Osaka University \\
Toyonaka, Osaka 560-0043, Japan}
\email{tomita@math.sci.osaka-u.ac.jp}



\begin{abstract}
The boundedness from 
$H^p \times L^2$ to $L^r$, $1/p+1/2=1/r$,
and from $H^p \times L^{\infty}$ to $L^p$ of 
bilinear pseudo-differential operators is proved 
under the assumption that their symbols are 
in the bilinear H\"ormander class 
$BS^m_{\rho,\rho}$, $0 \le \rho <1$, 
of critical order $m$, where $H^p$ is the Hardy space.  
This combined with the previous results of the same authors 
establishes the sharp boundedness 
from $H^p \times H^q$ to $L^r$, $1/p+1/q=1/r$, of those 
operators in the full range $0< p, q \le \infty$, 
where $L^r$ is replaced by $BMO$ if $r=\infty$. 
\end{abstract}

\maketitle

\section{Introduction}\label{section1}
This paper is a continuation of the paper \cite{Miyachi-Tomita-AIF}.  
We continue the study of the boundedness of bilinear pseudo-differential operators
with symbols in the so-called exotic classes. 
As for the background of this subject, see 
Introduction of \cite{Miyachi-Tomita-AIF}. 
Here we begin by recalling necessary definitions.

Let $m \in \R$ and $0 \le \delta \le \rho \le 1$.
We say that a function
$\sigma(x,\xi,\eta) \in C^{\infty}(\R^n \times \R^n \times \R^n)$
belongs to the bilinear H\"ormander symbol class
$BS^m_{\rho,\delta}=BS^m_{\rho,\delta}(\R^n)$
if for every triple of multi-indices
$\alpha,\beta,\gamma \in \N_0^n=\{0,1,2,\dots\}^n$
there exists a constant $C_{\alpha,\beta,\gamma}>0$ such that
\[
|\partial^{\alpha}_x\partial^{\beta}_{\xi}
\partial^{\gamma}_{\eta}\sigma(x,\xi,\eta)|
\le C_{\alpha,\beta,\gamma}
(1+|\xi|+|\eta|)^{m+\delta|\alpha|-\rho(|\beta|+|\gamma|)}.
\]
In this paper, we consider the class $BS^m_{\rho,\delta}$ 
with $0\le \rho=\delta <1$.

The bilinear pseudo-differential operator
$T_{\sigma}$, $\sigma \in BS^{m}_{\rho,\rho}$, is defined by
\[
T_{\sigma}(f,g)(x)
=\frac{1}{(2\pi)^{2n}}
\int_{\R^n \times \R^n}e^{i x \cdot(\xi+\eta)}
\sigma(x,\xi,\eta)\widehat{f}(\xi)\widehat{g}(\eta)\, d\xi d\eta,
\qquad f,g \in \calS(\R^n).
\]
If $X,Y,Z$ are function spaces on $\R^n$ 
equipped with quasi-norms 
$\|\cdot \|_{X},\,\|\cdot \|_{Y},\,\|\cdot \|_{Z}$ 
and if there exists a constant $A_{\sigma}$ such that 
the estimate 
$\|T_{\sigma}(f,g)\|_{Z}
\le A_{\sigma} \|f\|_{X} \|g\|_{Y}$
holds for all $f\in \calS \cap X$ and all $g\in \calS \cap Y$, 
then we shall simply say that 
$T_{\sigma}$ is bounded from $X\times Y$ to $Z$ and write 
$T_{\sigma}: X\times Y \to Z$. 
For the function spaces $X$ and $Y$, we consider 
the Hardy spaces $H^p$, $0<p\le \infty$. 
For $Z$, we consider the Lebesgue spaces 
$L^r$, $0<r<\infty$, or $BMO$. 
Notice that $H^p=L^p$ for $1<p\le \infty$. 
The definitions of $H^p$ and $BMO$ are given in Section 2.

For $0 \le \rho <1$
and for $0<p,q,r \le \infty$ satisfying $1/p+1/q=1/r$,
we write 
\begin{align*}
&
m_\rho(p,q)
=(1-\rho)m_0(p,q),
\\
&
m_0(p,q)
=-n \left(\max\left\{
\frac{1}{2}, \,
\frac{1}{p}, \,
\frac{1}{q}, \,
1-\frac{1}{r}, \, 
\frac{1}{r} -\frac{1}{2}
\right\}\right). 
\end{align*}
%
Observe that the region of $(1/p, 1/q)$, $0\le 1/p,1/q < \infty$, 
is divided into 5 regions, on each of which $m_{0}(p,q)$ is 
an affine function of $1/p$ and $1/q$ (see \cite[Introduction]{Miyachi-Tomita-AIF}).

The number $m_{\rho}(p,q)$ is the 
critical order as the following proposition shows.

\begin{referpropA}
Let $0 \le \rho <1$, $0<p,q,r \le \infty$, and suppose 
$1/p+1/q=1/r$. 
If $r<\infty$, then 
\[
m_{\rho}(p,q)
=
\sup \{m \in \R \,:\, 
T_{\sigma}: H^p \times H^q \to L^r 
\;\; \text{for all}\;\; \sigma \in BS^{m}_{\rho,\rho}
\}.  
\] 
When $p=q=r=\infty$, the above equality holds 
if we replace 
$H^p \times H^q \to L^r$ by 
$L^{\infty} \times L^{\infty} \to BMO $. 
\end{referpropA}

In fact, this proposition is a conclusion of several previous works:   
Michalowski-Rule-Staubach \cite{MRS} 
(for $(1/p, 1/q)$ in the triangle 
with vertices $(1/2, 1/2)$, $(1/2, 0)$, $(0, 1/2)$), 
B\'enyi-Bernicot-Maldonado-Naibo-Torres \cite{BBMNT} 
(in the range $1/p+1/q \le 1$), 
and Miyachi-Tomita \cite{Miyachi-Tomita-IUMJ, Miyachi-Tomita-AIF}  
(full range $0<p,q\le \infty$). 
For a proof of Proposition A, see \cite[Appendix A]{Miyachi-Tomita-AIF}.

It should be an interesting problem to prove the sharp boundedness, i.e., 
the boundedness $T_{\sigma}: H^p \times H^q \to L^r $, 
$1/r = 1/p+1/q$, 
with $L^r$ replaced by $BMO$ if $r=\infty$, 
for $\sigma \in BS^{m}_{\rho,\rho}$ with $m=m_{\rho}(p,q)$.

In the case $\rho=0$, this sharp boundedness was 
proved in \cite{Miyachi-Tomita-IUMJ}.

Recently, the authors proved the following theorem, which gives the 
sharp boundedness in the range $1\le p,q,r\le \infty$.

\begin{referthmB}[{\cite[Corollary 1.4]{Miyachi-Tomita-AIF}}]
Let $0 \le \rho <1$, $1 \le p,q,r \le \infty$, 
$1/p+1/q=1/r$, and $m=m_{\rho}(p,q)$. 
Then all bilinear pseudo-differential operators
with symbols in $BS^{m}_{\rho,\rho}(\R^n)$
are bounded from $L^p(\R^n) \times L^{q}(\R^n)$ to $L^r(\R^n)$, 
where $L^p(\R^n)$ and $L^q(\R^n)$ 
should be replaced by $H^1(\R^n)$ 
if $p=1$ or $q=1$ and 
$L^r(\R^n)$ should be replaced by $BMO(\R^n)$ if $r=\infty$.
\end{referthmB}

Here it should be mentioned that
Theorem B with $p=q=r=\infty$ 
and $0<\rho<1/2$ was also proved by  Naibo \cite{Naibo}.

Now the purpose of the present paper is to prove 
the sharp boundedness in the remaining cases 
and establish 
the sharp boundedness in the full range $0 <p, q, r \le \infty$.
The following is the conclusion of the present paper.

\begin{thm}\label{main-thm}
Let $0 \le \rho <1$, $0< p,q,r \le \infty$, 
$1/p+1/q=1/r$, and $m=m_{\rho}(p,q)$. 
Then all bilinear pseudo-differential operators
with symbols in $BS^{m}_{\rho,\rho}(\R^n)$
are bounded from $H^p(\R^n) \times H^{q}(\R^n)$ to $L^r(\R^n)$, 
where $L^r(\R^n)$ should be replaced by $BMO(\R^n)$ if $p=q=r=\infty$. 
\end{thm}

The above theorem follows with the aid of complex interpolation and symmetry 
if the sharp boundedness is proved in the following 5 cases: 
\begin{align*}
&
\mathrm{(i)}\ \ (p,q)=(\infty, \infty), 
\ m_{\rho}(\infty, \infty)= -(1-\rho)n; 
\\
&
\mathrm{(ii)}\ \ (p,q)=(2, 2), 
\ m_{\rho}(2, 2)= -(1-\rho)n/2; 
\\
&\mathrm{(iii)}\ \ (p,q)=(2, \infty), 
\ m_{\rho}(2, \infty)= -(1-\rho)n/2; 
\\
&\mathrm{(iv)}\ \ 0<p<1, \ q=2, 
\ m_{\rho}(p, 2)= -(1-\rho)n/p; 
\\
&
\mathrm{(v)}\ \ 0<p<1, \ q=\infty, 
\ m_{\rho}(p, \infty)= -(1-\rho)n/p.   
\end{align*}
(For the interpolation argument, see, e.g., 
\cite[Proof of Theorem 2.2]{BBMNT} 
or \cite[Proof of the `if' part of Theorem 1.1]{Miyachi-Tomita-IUMJ}.)  
By symbolic calculus of $BS^{m}_{\rho, \rho}$ as given 
by B\'enyi-Maldonado-Naibo-Torres \cite{BMNT} and by duality, 
the cases (ii) and (iii) are essentially the same (see, e.g., \cite[Section 5]{Miyachi-Tomita-AIF}).  
Thus Theorem \ref{main-thm} will follow if we prove 
(i), (ii)$=$(iii), (iv), and (v). 
Among these 4 critical cases, (i) and (ii)$=$(iii) are 
covered by Theorem B.  
Thus in order to prove Theorem \ref{main-thm} it is sufficient 
to prove (iv) and (v), which we shall state here as 
the following two theorems.

\begin{thm}\label{main-thm-1}
Let $0 \le \rho <1$, $0<p<1$, 
$1/p+1/2=1/r$, and $m=-(1-\rho)n/p$. 
Then all bilinear pseudo-differential operators
with symbols in $BS^{m}_{\rho,\rho}(\R^n)$
are bounded from $H^p(\R^n) \times L^{2}(\R^n)$ to $L^r(\R^n)$.
\end{thm}

\begin{thm}\label{main-thm-2}
Let $0 \le \rho <1$, $0<p<1$, and $m=-(1-\rho)n/p$. 
Then all bilinear pseudo-differential operators
with symbols in $BS^{m}_{\rho,\rho}(\R^n)$
are bounded from $H^p(\R^n) \times L^{\infty}(\R^n)$ to $L^p(\R^n)$.
\end{thm}

Here are some comments on the proofs of the theorems. 
First, although the case (i) 
(the sharp $L^{\infty} \times L^{\infty} \to BMO$ boundedness) 
was directly proved in \cite{Naibo, Miyachi-Tomita-AIF}, 
the argument of the present paper gives an alternate proof. 
In fact, by virtue of the symbolic calculus of $BS^{m}_{\rho, \rho}$ 
and by the duality $(H^1)^{\prime}=BMO$,  
the sharp $L^{\infty} \times L^{\infty} \to BMO$ boundedness 
is equivalent to the sharp 
$H^{1} \times L^{\infty} \to L^{1}$ boundedness 
and the latter follows from 
the cases (iii) and (v) (Theorem \ref{main-thm-2}) by interpolation. 
Secondly, for the proof of the case (iv) 
(Theorem \ref{main-thm-1}), 
the method of \cite{Miyachi-Tomita-IUMJ} given for 
$\rho=0$ does not seem to work for $0<\rho<1$. 
Our proof of Theorem \ref{main-thm-1} is based on a new method, 
which covers $\rho=0$ and $0 <\rho <1$ simultaneously.  
Finally, the case (v) (Theorem \ref{main-thm-2}) is rather easy. 
In fact, by freezing $g$ of $T_{\sigma}(f,g)$ 
we can follow the argument used in the case of 
linear pseudo-differential operators.

The contents of this paper are as follows.
In Section \ref{section2},
we recall some preliminary facts.
In Sections \ref{section3} and \ref{section4},
we prove Theorems \ref{main-thm-1} and \ref{main-thm-2},
respectively.

\section{Preliminaries}\label{section2}

For two nonnegative quantities $A$ and $B$,
the notation $A \lesssim B$ means that
$A \le CB$ for some unspecified constant $C>0$,
and $A \approx B$ means that
$A \lesssim B$ and $B \lesssim A$.
We denote by $\ichi_S$ the characteristic function of a set $S$,
and by $|S|$ the Lebesgue measure of a measurable set $S$ in $\R^n$.

Let $\calS(\R^n)$ and $\calS'(\R^n)$ be the Schwartz space of
rapidly decreasing smooth functions and the space of
tempered distributions, respectively.
We define the Fourier transform $\calF f$
and the inverse Fourier transform $\calF^{-1}f$
of $f \in \calS(\R^n)$ by
\[
\calF f(\xi)
=\widehat{f}(\xi)
=\int_{\R^n}e^{-ix\cdot\xi} f(x)\, dx
\quad \text{and} \quad
\calF^{-1}f(x)
=\frac{1}{(2\pi)^n}
\int_{\R^n}e^{ix\cdot \xi} f(\xi)\, d\xi.
\]
For $m \in L^{\infty}(\R^n)$,
the linear Fourier multiplier operator $m(D)$
is defined by
\[
m(D)f(x)
=\calF^{-1}[m\widehat{f}](x)
=\frac{1}{(2\pi)^n}\int_{\R^n}
e^{ix\cdot\xi}m(\xi)\widehat{f}(\xi)\, d\xi,
\quad
f \in \calS(\R^n).
\]

We recall the definitions and some properties of 
Hardy spaces and the space $BMO$ on $\R^n$
(see, e.g., \cite[Chapters 3 and 4]{Stein}).
Let $0 < p \le \infty$, and let $\phi \in \calS(\R^n)$ be such that
$\int_{\R^n}\phi(x)\, dx \neq 0$. 
Then the Hardy space $H^p(\R^n)$ consists of all $f \in \calS'(\R^n)$
such that
\[
\|f\|_{H^p}=\|\sup_{0<t<\infty}|\phi_t*f|\|_{L^p}<\infty,
\]
where $\phi_t(x)=t^{-n}\phi(x/t)$.
It is known that $H^p(\R^n)$ does not depend 
on the choice of the function $\phi$
and $H^p(\R^n)=L^p(\R^n)$ for $1<p \le \infty$.
For $0 <p \le 1$, a function $a$ on $\R^n$ 
is called an $H^p$-atom 
if there exists a cube $Q=Q_a$ such that
\begin{equation}\label{atomQ}
\mathrm{supp}\, a \subset Q,
\qquad \|a\|_{L^{\infty}} \le |Q|^{-1/p},
\qquad \int_{\R^n}x^{\alpha}\, a(x)\, dx=0,
\ |\alpha| \le L-1,
\end{equation}
where 
$L$ is any fixed integer
satisfying $L > n/p-n$ 
(\cite[p.112]{Stein}).
It is known that every $f \in H^p(\R^n)$ can be written
as
\[
f=\sum_{i=1}^{\infty}\lambda_i a_i
\quad \text{in} \quad \calS'(\R^n),
\]
where $\{a_i\}$ is a collection of $H^p$-atoms
and $\{\lambda_i\}$ is a sequence of complex numbers with
$\sum_{i=1}^{\infty}|\lambda_i|^p<\infty$.
Moreover,
\[
\|f\|_{H^p}
\approx
\inf \left(\sum_{i=1}^{\infty}|\lambda_i|^p\right)^{1/p},
\]
where the infimum is taken over all representations of $f$.
The space $BMO(\R^n)$ consists of
all locally integrable functions $f$ on $\R^n$ such that
\[
\|f\|_{BMO}
=\sup_{Q}\frac{1}{|Q|}
\int_{Q}|f(x)-f_Q|\, dx<\infty, 
\]
where $f_Q$ is the average of $f$ on $Q$
and the supremum is taken over all cubes $Q$ in $\R^n$.
It is known that the dual space of
$H^1(\R^n)$ is $BMO(\R^n)$.

\section{Proof of Theorem \ref{main-thm-1}}\label{section3}

In this section,
we shall prove Theorem \ref{main-thm-1}. 
We assume $0 \le \rho<1$, $0<p<1$, $1/p+1/2=1/r$,
$m=-(1-\rho)n/p$, and $\sigma \in BS^{m}_{\rho,\rho}$, 
and prove the $H^p \times L^2 \to L^r$ boundedness of $T_{\sigma}$.

We first observe that the desired boundedness follows
if we prove the following:
for an $H^p$-atom $a$ and an $L^2$-function $g$
there exist a function $\widetilde{a}$ depending only on $a$
and a function $\widetilde{g}$ depending only on $g$ such that
\begin{equation}\label{thm1-goal1}
|T_{\sigma}(a,g)(x)|
\lesssim
\widetilde{a}(x)\widetilde{g}(x),
\quad
\|\widetilde{a}\|_{L^p} \lesssim 1,
\quad
\|\widetilde{g}\|_{L^2} \lesssim \|g\|_{L^2}.
\end{equation}
In fact,
if this is proved,
we can deduce the $H^p \times L^2 \to L^r$
boundedness of $T_{\sigma}$ as follows.
Given $f \in H^p$,
we decompose it as
\[
f=\sum_{i}\lambda_i a_i,
\qquad
\left(\sum_{i}|\lambda_i|^p\right)^{1/p}
\lesssim \|f\|_{H^p},
\]
where $a_i$, $i \ge 1$, are $H^p$-atoms.
Then, taking the functions $\widetilde{a}_i$ and $\widetilde{g}$
satisfying \eqref{thm1-goal1} for $a=a_i$,
we have
\begin{align*}
\|T_{\sigma}(f,g)\|_{L^r}
&=\left\|\sum_{i}\lambda_i T_{\sigma}(a_i,g)\right\|_{L^r}
\lesssim
\left\|\left(\sum_{i}|\lambda_i|\widetilde{a}_i\right)
\widetilde{g}\right\|_{L^r}
\\
&\le \left\|\sum_{i}|\lambda_i|\widetilde{a}_i\right\|_{L^p}
\|\widetilde{g}\|_{L^2}
\lesssim
\left(\sum_{i}|\lambda_i|^p\right)^{1/p}\|g\|_{L^2}
\lesssim \|f\|_{H^p}\|g\|_{L^2}.
\end{align*}
(This argument was already used in \cite{Miyachi-Tomita-RMI}. 
The idea goes back to \cite{G-K}.)

Let $a$ be an $H^p$-atom satisfying \eqref{atomQ} with $L>n/p-n$.
We denote by $c_Q$ the center of $Q$,
by $\ell(Q)$ the side length of $Q$,
and by $Q^*$ the cube with the same center as $Q$
but expanded by a factor of $2\sqrt{n}$.
To obtain \eqref{thm1-goal1},
we shall prove
\begin{alignat}{3}
&|T_{\sigma}(a,g)(x)|\ichi_{(Q^*)^c}(x)
\lesssim u(x)v(x),
&\quad
&\|u\|_{L^p} \lesssim 1,
&\quad
&\|v\|_{L^2} \lesssim \|g\|_{L^2},
\label{thm1-goal1-outside}
\\
&|T_{\sigma}(a,g)(x)|\ichi_{Q^*}(x)
\lesssim u'(x)v'(x),
&\quad
&\|u'\|_{L^p} \lesssim 1,
&\quad
&\|v'\|_{L^2} \lesssim \|g\|_{L^2},
\label{thm1-goal1-inside}
\end{alignat}
where $u, u'$ depend only on $a$
and $v,v'$ depend only on $g$.
Once \eqref{thm1-goal1-outside} and \eqref{thm1-goal1-inside}
are proved, we can take $u+u'$ and $v+v'$ as $\widetilde{a}$
and $\widetilde{g}$ in \eqref{thm1-goal1}.

Let $\psi_0 \in \calS(\R^d)$ be such that
$\psi_0=1$ on $\{\zeta \in \R^d \,:\, |\zeta| \le 1\}$
and $\supp \psi_0 \subset \{\zeta \in \R^d \,:\, |\zeta| \le 2\}$,
and set $\psi(\zeta)=\psi_0(\zeta)-\psi_0(2\zeta)$
and $\psi_j(\zeta)=\psi(\zeta/2^j)$,
$j \ge 1$.
Then
\begin{equation}\label{dyadic-dec}
\supp \psi_j \subset \{\zeta \in \R^d \,:\, 
2^{j-1} \le |\zeta| \le 2^{j+1}\},
\ j \ge 1, \qquad
\sum_{j=0}^{\infty}\psi_j(\zeta)=1,
\ \zeta \in \R^d.
\end{equation}
We also use functions
$\widetilde{\psi}_0, \widetilde{\psi} \in \calS(\R^n)$
satisfying
$\widetilde{\psi}_0=1$ on $\{\eta \in \R^n \,:\, |\eta| \le 4\}$,
$\supp \widetilde{\psi}_0 \subset \{\eta \in \R^n \,:\, |\eta| \le 8\}$,
$\widetilde{\psi}=1$ on $\{\eta \in \R^n \,:\, 1/4 \le |\eta| \le 4\}$,
and $\supp \widetilde{\psi} \subset
\{\eta \in \R^n \,:\, 1/8 \le |\eta| \le 8\}$,
and set $\widetilde{\psi}_\ell(\eta)=\widetilde{\psi}(\eta/2^{\ell})$,
$\ell \ge 1$.
In order to obtain \eqref{thm1-goal1-outside}
and \eqref{thm1-goal1-inside},
we decompose $T_{\sigma}(a,g)$ as
\begin{equation}\label{goal1-split-sum}
T_{\sigma}(a,g)(x)
=\sum_{j=0}^{\infty}\sum_{\ell=0}^{\infty}T_{\sigma_{j,\ell}}(a,g)(x)
=\sum_{\substack{j, \ell \ge 0 \\ \ell \le j(1-\rho)+2}}
T_{\sigma_{j,\ell}}(a,g_{j,\ell})(x)
\end{equation}
with
\[
\sigma_{j,\ell}(x,\xi,\eta)=\sigma(x,\xi,\eta)
\Psi_{j}(\xi,\eta)\psi_{\ell}(\eta/2^{j\rho})
\]
and
\[
g_{j,\ell}(x)=\widetilde{\psi}_{\ell}(D/2^{[j\rho]})g(x),
\]
where
$\Psi_j$ and $\psi_\ell$
are as in \eqref{dyadic-dec} with $d=2n$ and $d=n$ respectively,
and $[j\rho]$ is the integer part of $j\rho$.
Here, we used the fact 
\[
\sum_{j \ge 0}\Psi_j(\xi,\eta)=\sum_{j \ge 0}\sum_{\ell \ge 0}
\Psi_j(\xi,\eta)\psi_{\ell}(\eta/2^{j\rho})=1,
\quad (\xi,\eta) \in \R^n \times \R^n,
\]
in the first equality of \eqref{goal1-split-sum}, and the facts 
\[
\Psi_j(\xi,\eta)\psi_{\ell}(\eta/2^{j\rho})=0,
\quad \ell>j(1-\rho)+2,
\]
and
\[
\psi_{\ell}(\eta/2^{j\rho})=\psi_{\ell}(\eta/2^{j\rho})
\widetilde{\psi}_{\ell}(\eta/2^{[j\rho]}),
\quad \ell \ge 0,
\]
in the second equality of \eqref{goal1-split-sum}.
We write the partial inverse Fourier transform of 
$\sigma_{j,\ell}(x,\xi, \eta)$ with respect to $(\xi, \eta)$ 
as
\[
K_{j,\ell}(x,y,z)
=\frac{1}{(2\pi)^{2n}}
\int_{(\R^n)^2}e^{i(y\cdot\xi+z\cdot\eta)}
\sigma_{j,\ell}(x,\xi,\eta)\, d\xi d\eta,
\quad
x,y,z \in \R^n,
\]
and then
\[
T_{\sigma_{j,\ell}}(a,g_{j,\ell})(x)
=\int_{(\R^n)^2}
K_{j,\ell}(x,x-y,x-z)a(y)g_{j,\ell}(z)\, dydz.
\]

Notice that $\sigma_{j,\ell}$ satisfies the following: 
\begin{align}
& \supp \sigma_{j,\ell} (x, \cdot, \cdot) 
\subset 
\{ |\xi|\le 2^{j+1}, \; |\eta|\le 2^{j\rho + \ell +1}\}, 
\label{sigmajell-1}
\\
& 
1+ |\xi|+ |\eta| \approx 2^{j} \quad \text{on}\quad 
\supp \sigma_{j,\ell} (x, \cdot, \cdot), 
\label{sigmajell-2}
\\
&
|\partial_{x}^{\alpha}\partial_{\xi}^{\beta}\partial_{\eta}^{\gamma}
\sigma_{j,\ell} (x, \xi, \eta) |
\lesssim 
2^{j m}\, 2^{j \rho (|\alpha|- |\beta| - |\gamma|)}. 
\label{sigmajell-3}
\end{align}

\begin{proof}[Proof of \eqref{thm1-goal1-outside}]
Let $x \not\in Q^*$.
Using the moment condition on $a$ and Taylor's formula, we have
\begin{equation}\label{goal1-moment}
\begin{split}
&T_{\sigma_{j,\ell}}(a,g_{j,\ell})(x)
=\int_{(\R^n)^2}
K_{j,\ell}(x,x-y,x-z)a(y)g_{j,\ell}(z)\, dydz
\\
&=\int_{(\R^n)^2}
\left(K_{j,\ell}(x,x-y,x-z)-\sum_{|\alpha|<L}
\frac{(c_Q-y)^{\alpha}}{\alpha!}
K_{j,\ell}^{(\alpha,0)}(x,x-c_Q,x-z)\right)
\\
&\qquad \qquad \times a(y)g_{j,\ell}(z)\, dydz
\\
&=L\sum_{|\alpha|=L}
\int_{\substack{y\in Q \\ z \in \R^n \\ 0<t<1}}
(1-t)^{L-1}\frac{(c_Q-y)^{\alpha}}{\alpha!}
K_{j,\ell}^{(\alpha,0)}(x,x-[c_Q,y]_t,x-z)
\\
&\qquad \qquad \times  a(y)g_{j,\ell}(z)\, dydzdt,
\end{split}
\end{equation}
where
$K_{j,\ell}^{(\alpha,0)}(x,y,z)=\partial_y^{\alpha}K_{j,\ell}(x,y,z)$
and $[c_Q,y]_t=c_Q+t(y-c_Q)$.
It follows from the size condition on $a$ that
\[
|T_{\sigma_{j,\ell}}(a,g_{j,\ell})(x)|
\lesssim
\ell(Q)^{L-n/p}
\sum_{|\alpha|=L}
\int_{\substack{y\in Q \\ z \in \R^n \\ 0<t<1}}
|K_{j,\ell}^{(\alpha,0)}(x,x-[c_Q,y]_t,x-z)
g_{j,\ell}(z)|\, dydzdt.
\]
Let $M$ and $M'$ be integers satisfying $M>n/p-n/2$ and $M'>n/2$.
Since $|x-c_Q| \approx |x-[c_Q,y]_t|$
for $x \not\in Q^*$, $y \in Q$ and $0<t<1$,
Schwarz's inequality with respect to the $z$-variable gives 
\begin{align*}
&(1+2^{j\rho}|x-c_Q|)^{M}
|T_{\sigma_{j,\ell}}(a,g_{j,\ell})(x)|
\\
&\lesssim
\ell(Q)^{L-n/p}
\sum_{|\alpha|=L}
\int_{\substack{y\in Q \\ z \in \R^n \\ 0<t<1}}
\\
&\qquad \times (1+2^{j\rho}|x-[c_Q,y]_t|)^{M}
|K_{j,\ell}^{(\alpha,0)}(x,x-[c_Q,y]_t,x-z)
g_{j,\ell}(z)|\, dydzdt
\\
&\le \ell(Q)^{L-n/p}
\sum_{|\alpha|=L}
\int_{\substack{y\in Q \\ 0<t<1}}
\\
&\qquad  \times
\left\|(1+2^{j\rho}|x-[c_Q,y]_t|)^{M}(1+2^{j\rho}|x-z|)^{M'}
K_{j,\ell}^{(\alpha,0)}(x,x-[c_Q,y]_t,x-z)\right\|_{L^2_z}
\\
&\qquad \times
\left\|(1+2^{j\rho}|x-z|)^{-M'} g_{j,\ell}(z)\right\|_{L^2_z}\, dydt.
\end{align*}
Thus, by writing
\begin{align*}
h_{j,\ell}^{(Q,L)}(x)
&=2^{-j\rho n/2}\ell(Q)^{L-n/p}
\sum_{|\alpha|=L}\sum_{|\beta| \le M}\sum_{|\gamma| \le M'}
\int_{\substack{y\in Q \\ 0<t<1}}
\\
&\qquad  \times
\left\|(2^{j\rho}(x-[c_Q,y]_t))^{\beta}(2^{j\rho}z)^{\gamma}
K_{j,\ell}^{(\alpha,0)}(x,x-[c_Q,y]_t,z)\right\|_{L^2_z}\, dydt
\end{align*}
and
\begin{equation}\label{tildegjell}
\widetilde{g}_{j,\ell}(x)
=2^{j\rho n/2}\left\|(1+2^{j\rho}|x-\cdot|)^{-M'}
g_{j,\ell}(\cdot)\right\|_{L^2}, 
\end{equation}
we have
\begin{equation}\label{goal1-estimate-moment}
|T_{\sigma_{j,\ell}}(a,g_{j,\ell})(x)|
\lesssim
(1+2^{j\rho}|x-c_Q|)^{-M}
h_{j,\ell}^{(Q,L)}(x)
\widetilde{g}_{j,\ell}(x).
\end{equation}

We shall estimate the $L^2$-norm of $h_{j,\ell}^{(Q,L)}$.
By Minkowski's inequality for integrals,
\begin{align*}
\|h_{j,\ell}^{(Q,L)}\|_{L^2}
&\le
2^{-j\rho n/2}\ell(Q)^{L-n/p}
\sum_{|\alpha|=L}\sum_{|\beta| \le M}\sum_{|\gamma| \le M'}
\int_{\substack{y\in Q \\ 0<t<1}}
\\
&\quad  \times
\left\|(2^{j\rho}(x-[c_Q,y]_t))^{\beta}(2^{j\rho}z)^{\gamma}
K_{j,\ell}^{(\alpha,0)}(x,x-[c_Q,y]_t,z)\right\|_{L^2_{x,z}}\, dydt.
\end{align*}
The function in the above $\|\cdot\|_{L^2_{x,z}}$
can be written as
\begin{align*}
&(2^{j\rho}(x-[c_Q,y]_t))^{\beta}(2^{j\rho}z)^{\gamma}
K_{j,\ell}^{(\alpha,0)}(x,x-[c_Q,y]_t,z)
\\
&=C_{\alpha,\beta,\gamma}
\int_{(\R^n)^2}
e^{i\{(x-[c_Q,y]_t)\cdot\xi+z\cdot\eta\}}
(2^{j\rho}\partial_{\xi})^{\beta}(2^{j\rho}\partial_{\eta})^{\gamma}
\left[\xi^{\alpha}\sigma_{j,\ell}(x,\xi,\eta)\right] d\xi d\eta
\\
&=C_{\alpha,\beta,\gamma}
\int_{(\R^n)^2}
e^{i\{x\cdot\xi+z\cdot\eta\}}
(2^{j\rho}\partial_{\xi})^{\beta}(2^{j\rho}\partial_{\eta})^{\gamma}
\left[\xi^{\alpha}\sigma_{j,\ell}(x,\xi,\eta)\right]
\\
&\qquad \qquad \qquad \times
e^{-i[c_Q,y]_t \cdot \xi}\psi_0(\xi/2^{j+1})\psi_0(\eta/2^{j\rho+\ell+1})
\, d\xi d\eta,
\end{align*}
where $\psi_0$ is as in \eqref{dyadic-dec} with $d=n$
and we used \eqref{sigmajell-1}. 
From \eqref{sigmajell-2} and \eqref{sigmajell-3}, 
we see that
\[
\left|\partial_x^{\alpha'}\partial_{\xi}^{\beta'}\partial_{\eta}^{\gamma'}
(2^{j\rho}\partial_{\xi})^{\beta}(2^{j\rho}\partial_{\eta})^{\gamma}
\left[\xi^{\alpha}\sigma_{j,\ell}(x,\xi,\eta)\right]\right|
\lesssim 2^{j(|\alpha|+m)}
(1+|\xi|+|\eta|)^{\rho(|\alpha'|-|\beta'|-|\gamma'|)}.
\]
Hence, the Calder\'on-Vaillancourt theorem on $\R^{2n}$ (\cite{CV}) 
and Plancherel's theorem give
\begin{equation*}
\begin{split}
&\left\|(2^{j\rho}(x-[c_Q,y]_t))^{\beta}(2^{j\rho}z)^{\gamma}
K_{j,\ell}^{(\alpha,0)}(x,x-[c_Q,y]_t,z)\right\|_{L^2_{x, z}}
\\
&\lesssim
2^{j(|\alpha| +m )}
\left\|e^{-i[c_Q,y]_t \cdot \xi}\psi_0(\xi/2^{j+1})
\psi_0(\eta/2^{j\rho+\ell+1})
\right\|_{L^2_{\xi,\eta}}
\\
&\approx  2^{j(|\alpha|+m + n/2)}
2^{(j \rho + \ell)n/2},
\end{split}
\end{equation*}
which implies
\begin{equation}\label{goal1-estimate-h_j^{(Q,L)}}
\begin{split}
\|h_{j,\ell}^{(Q,L)}\|_{L^2}
&\lesssim
2^{-j\rho n/2}\ell(Q)^{L-n/p}
\sum_{|\alpha|=L}
2^{j(|\alpha|+m + n/2)}
2^{(j \rho + \ell)n/2}
|Q|
\\
&\approx \left(2^j \ell(Q)\right)^{L-n/p+n}2^{j\rho n(1/p-1/2)}
2^{(\ell-j(1-\rho))n/2}.
\end{split}
\end{equation}

If we do not use the moment condition on $a$ in \eqref{goal1-moment},
a similar argument yields
\begin{equation}\label{goal1-estimate-not-moment}
|T_{\sigma_{j,\ell}}(a,g_{j,\ell})(x)|
\lesssim
(1+2^{j\rho}|x-c_Q|)^{-M}
h_{j,\ell}^{(Q,0)}(x)
\widetilde{g}_{j,\ell}(x)
\end{equation}
with
\begin{align*}
h_{j,\ell}^{(Q,0)}(x)
&=2^{-j\rho n/2}\ell(Q)^{-n/p}
\sum_{|\beta| \le M}\sum_{|\gamma| \le M'}
\int_{y\in Q}
\\
&\qquad \times
\left\|(2^{j\rho}(x-y))^{\beta}(2^{j\rho}z)^{\gamma}
K_{j,\ell}(x,x-y,z)\right\|_{L^2_z}\, dy
\end{align*}
and
\begin{equation}\label{goal1-estimate-h_j^{(Q,0)}}
\|h_{j,\ell}^{(Q,0)}\|_{L^2}
\lesssim
\left(2^j \ell(Q)\right)^{-n/p+n}2^{j\rho n(1/p-1/2)}
2^{(\ell-j(1-\rho))n/2}.
\end{equation}

Combining \eqref{goal1-estimate-moment}
and \eqref{goal1-estimate-not-moment}, we have
\[
|T_{\sigma_{j,\ell}}(a,g_{j,\ell})(x)|
\lesssim
u_{j,\ell}(x)\widetilde{g}_{j,\ell}(x)
\]
with
\[
u_{j,\ell}(x)=
(1+2^{j\rho}|x-c_Q|)^{-M}
\min\left\{h_{j,\ell}^{(Q,L)}(x), h_{j,\ell}^{(Q,0)}(x)\right\}. 
\]
We take an $\epsilon$ satisfying $0<\epsilon < n/2$ and set
\begin{equation*}
u(x)=\left(\sum_{\ell \le j(1-\rho)+2}
2^{-(\ell -j (1-\rho))2 \epsilon }\, 
u_{j,\ell}(x)^2\right)^{1/2} 
\end{equation*}
and 
\begin{equation}\label{goal1-uv2}
v(x)=\left(\sum_{\ell \le j(1-\rho)+2}
2^{(\ell -j (1-\rho))2 \epsilon }\, 
\widetilde{g}_{j,\ell}(x)^2\right)^{1/2}. 
\end{equation}
(The number $\epsilon$ can be chosen arbitrarily in the range $0<\epsilon <n/2$. 
For example $\epsilon = n/4$ suffices.)  
Then Schwarz's inequality gives
\[
|T_{\sigma}(a,g)(x)|
\le \sum_{\ell \le j(1-\rho)+2}
\left|T_{\sigma_{j,\ell}}(a,g_{j,\ell})(x)\right|
\lesssim
\sum_{\ell \le j(1-\rho)+2}u_{j,\ell}(x)\widetilde{g}_{j,\ell}(x)
\le u(x)v(x)
\]
for $x \not\in Q^*$.
Certainly the function $u$ depends only on $a$ and the 
function $v$ depends only on $g$. 
In the rest of the argument, we shall prove that
$\|u\|_{L^p} \lesssim 1$ and $\|v\|_{L^2} \lesssim \|g\|_{L^2}$,
which will complete the proof of \eqref{thm1-goal1-outside}.

First we shall prove $\|u\|_{L^p} \lesssim 1$. 
By H\"older's inequality with $1/p=1/q+1/2$ and by 
\eqref{goal1-estimate-h_j^{(Q,L)}}
and \eqref{goal1-estimate-h_j^{(Q,0)}},
\begin{align*}
\|u_{j,\ell}\|_{L^p}
&\le 
\|(1+2^{j\rho}|\cdot-c_Q|)^{-M}\|_{L^q}
\left\|\min\left\{h_{j,\ell}^{(Q,L)}, h_{j,\ell}^{(Q,0)}\right\}\right\|_{L^2}
\\
&\lesssim
2^{-j\rho n/q}\min\left\{
\|h_{j,\ell}^{(Q,L)}\|_{L^2}, \|h_{j,\ell}^{(Q,0)}\|_{L^2}\right\}
\\
&\lesssim 
2^{(\ell - j (1-\rho))n/2} 
\min\left\{\left(2^j \ell(Q)\right)^{L-n/p+n},
\left(2^j \ell(Q)\right)^{-n/p+n}\right\}.
\end{align*}
Thus
\begin{align*}
&
\|u\|_{L^p}^{p}
=
\left\|\left(\sum_{\ell \le j(1-\rho)+2}
2^{-(\ell - j (1-\rho)) 2 \epsilon }\, 
u_{j,\ell}^2\right)^{1/2}
\right\|_{L^p}^p
\\
&\le
\sum_{\ell \le j(1-\rho)+2}
\left(2^{-(\ell - j (1-\rho)) \epsilon }\|u_{j,\ell}\|_{L^p}\right)^p
\\
&\lesssim
\sum_{j=0}^{\infty}
\left(\min\left\{\left(2^j \ell(Q)\right)^{L-n/p+n},
\left(2^j \ell(Q)\right)^{-n/p+n}\right\}\right)^p
\left(\sum_{\ell=0}^{[j(1-\rho)]+2}
2^{(\ell-j(1-\rho))(n/2 - \epsilon )p}\right)
\\
&\approx  
\sum_{j=0}^{\infty}
\left(\min\left\{\left(2^j \ell(Q)\right)^{L-n/p+n},
\left(2^j \ell(Q)\right)^{-n/p+n}\right\}\right)^p
\lesssim 1,
\end{align*}
where the last $\lesssim$ holds because $L-n/p+n>0$ and $-n/p+n<0$.

Next, to prove $\|v\|_{L^2} \lesssim \|g\|_{L^2}$, 
observe that 
$\left\|\widetilde{g}_{j,\ell}\right\|_{L^2}\approx 
\left\|g_{j,\ell}\right\|_{L^2}$. 
Hence 
\begin{equation*}
\|v\|_{L^2}^2
=\sum_{\ell \le j(1-\rho) +2 } 
2^{(\ell-j(1-\rho))2 \epsilon }\, 
\left\|\widetilde{g}_{j,\ell}\right\|_{L^2}^2 
\approx 
\sum_{j=0}^{\infty}
\sum_{\ell=0}^{[j(1-\rho)]+2}
2^{(\ell-j(1-\rho))2 \epsilon }\, 
\left\|g_{j,\ell}\right\|_{L^2}^2. 
\end{equation*}
We divide the sum over $\ell$
into two parts $\ell=0$ and $\ell \ge 1$.  
For the terms with $\ell=0$, Young's inequality 
gives $\left\|g_{j,0}\right\|_{L^2}\le 
\|\calF^{-1}\widetilde{\psi}_0\|_{L^1}\|g\|_{L^2}
\approx \|g\|_{L^2}$ and thus 
\begin{equation*}
\sum_{j=0}^{\infty}
2^{-j(1-\rho)2 \epsilon }
\left\|g_{j,0}\right\|_{L^2}^2 
\lesssim 
\sum_{j=0}^{\infty}
2^{-j(1-\rho)2 \epsilon }
\|g\|_{L^2}^2 
\approx 
\|g\|_{L^2}^2 
\end{equation*}
since $\rho <1$. 
For the terms with $\ell \ge 1$, we have 
$g_{j,\ell}=\widetilde{\psi}(D/ 2^{[j\rho]+ \ell}) g$ 
and hence, 
by a change of variables, we have 
\begin{align*}
&
\sum_{j=0}^{\infty}
\sum_{\ell =1}^{[j(1- \rho)]+2}
2^{(\ell -j(1-\rho)) 2 \epsilon }
\left\|g_{j,\ell}\right\|_{L^2}^2 
\\
&= 
\sum_{j=0}^{\infty}
\sum_{k=[j\rho]+1}^{[j\rho ]+ [j(1-\rho)]+2}
2^{(k-[j\rho]-j(1-\rho))2 \epsilon }
\left\|\widetilde{\psi}(D/ 2^{k}) g\right\|_{L^2}^2 
\\
&\lesssim
\sum_{j=0}^{\infty}
\sum_{k=1}^{j+2}
2^{(k -j)2 \epsilon }
\left\|\widetilde{\psi}(D/2^{k})g\right\|_{L^2}^2
=\sum_{k=1}^{\infty}\sum_{j=\max\{0,\, k-2\}}^{\infty}
2^{(k -j)2 \epsilon }
\left\|\widetilde{\psi}(D/2^{k})g\right\|_{L^2}^2
\\
&\approx 
\sum_{k=1}^{\infty}
\left\|\widetilde{\psi}(D/2^{k})g\right\|_{L^2}^2
\lesssim \|g\|_{L^2}^2,
\end{align*}
where the last $\lesssim $ follows from the fact
$\mathrm{supp}\, \widetilde{\psi}$ is included in an annulus. 
Therefore, we obtain $\|v\|_{L^2} \lesssim \|g\|_{L^2}$.
\end{proof}

\begin{proof}[Proof of \eqref{thm1-goal1-inside}]
Take an $M^{\prime} > n/2$. 
By Schwarz's inequality,
\begin{align*}
&|T_{\sigma_{j,\ell}}(a,g_{j,\ell})(x)|
\le |Q|^{-1/p}\int_{(\R^n)^2}
|K_{j,\ell}(x,x-y,x-z)g_{j,\ell}(z)|\, dydz
\\
&\le
|Q|^{-1/p}
\left\|(1+2^{j\rho}|x-y|)^{M^{\prime}}(1+2^{j\rho}|x-z|)^{M^{\prime}}
K_{j,\ell}(x,x-y,x-z)\right\|_{L^2_{y,z}}
\\
&\qquad \times
\left\|(1+2^{j\rho}|x-y|)^{-M^{\prime}}(1+2^{j\rho}|x-z|)^{-M^{\prime}}
g_{j,\ell}(z)\right\|_{L^2_{y,z}}
\end{align*}
For the first $L^2_{y,z}$ norm above, 
we use Plancherel's theorem, \eqref{sigmajell-3} and 
\eqref{sigmajell-1} to obtain 
\begin{align*}
&\left\|(1+2^{j\rho}|y|)^{M^{\prime}}(1+2^{j\rho}|z|)^{M^{\prime}}
K_{j,\ell}(x,y,z)\right\|_{L^2_{y,z}} 
\\
&\approx 
\sum_{|\beta|\le M^{\prime}} 
\sum_{|\gamma|\le M^{\prime}} 
\left\|
(2^{j\rho}\partial_{\xi})^{\beta}
(2^{j\rho}\partial_{\eta})^{\gamma}
\sigma_{j,\ell} (x, \xi, \eta)
\right\|_{L^2_{\xi, \eta}} 
\\
&\lesssim 
2^{j (m+n/2)} 2^{(j\rho + \ell) n/2} 
\end{align*}
for all $x \in \R ^{n}$. 
As for the second $L^2_{y,z}$ norm, we have 
\begin{align*}
&
\left\|(1+2^{j\rho}|x-y|)^{-M^{\prime}}
(1+2^{j\rho}|x-z|)^{-M^{\prime}}
g_{j,\ell}(z)\right\|_{L^2_{y,z}}
\\
&\approx 
2^{-j \rho n/2}
\left\|(1+2^{j\rho}|x-z|)^{-M'}
g_{j,\ell}(z)\right\|_{L^2_{z}}
=2^{-j \rho n}\, \widetilde{g}_{j, \ell} (x),  
\end{align*}
where $\widetilde{g}_{j, \ell}$ is defined by \eqref{tildegjell}. 
Thus 
\begin{equation}\label{aaaa}
\begin{split}
|T_{\sigma_{j,\ell}}(a,g_{j,\ell})(x)|
&\lesssim
|Q|^{-1/p}
2^{j (m+n/2)} 2^{(j \rho + \ell ) n/2} 
2^{-j \rho n}\, 
\widetilde{g}_{j,\ell} (x)
\\
&
=
|Q|^{-1/p}
2^{-j (1-\rho)n (1/p-1))} 
2^{(\ell-j(1-\rho))n/2}\, 
\widetilde{g}_{j,\ell} (x). 
\end{split}
\end{equation}
Since 
\begin{align*}
&\sum_{\ell \le j(1-\rho)+2}
\left(
2^{-j(1-\rho)n(1/p-1)}
2^{(\ell-j(1-\rho)) (n/2 - \epsilon )}\right)^2
\\
&=\sum_{j=0}^{\infty}2^{-2j(1-\rho)n(1/p-1)}
\left(\sum_{\ell=0}^{[j(1-\rho)]+2}2^{(\ell-j(1-\rho))(n - 2 \epsilon )}\right)
\\
&\approx 
\sum_{j=0}^{\infty}2^{-2j(1-\rho)n(1/p-1)} 
\approx 1,  
\end{align*}
the estimate \eqref{aaaa} together with 
Schwarz's inequality gives 
\[
|T_{\sigma}(a,g)(x)|
\le \sum_{\ell \le j(1-\rho)+2} 
|T_{\sigma_{j,\ell}}(a,g_{j,\ell})(x)| 
\lesssim |Q|^{-1/p}v(x), 
\]
where $v$ is the function defined by \eqref{goal1-uv2}. 
In particular 
\[|T_{\sigma}(a,g)(x)| \ichi_{Q^*}(x) \lesssim |Q|^{-1/p}\ichi_{Q^*}(x) v(x). 
\]  
In Proof of \eqref{thm1-goal1-outside}, we have proved 
that $\|v\|_{L^2}\lesssim \|g\|_{L^2}$. 
Therefore, we can take $|Q|^{-1/p}\ichi_{Q^*}$ and $v$
as $u'$ and $v'$ in \eqref{thm1-goal1-inside}.
\end{proof}

\section{Proof of Theorem \ref{main-thm-2}}\label{section4}

In this section,
we shall prove Theorem \ref{main-thm-2},
namely the $H^p \times L^\infty \to L^p$ boundedness of $T_{\sigma}$,
where $0 \le \rho<1$, $0<p<1$,
$m=-(1-\rho)n/p$, and $\sigma \in BS^{m}_{\rho,\rho}$.

By the usual argument using the atomic decomposition for $H^p$,
the desired boundedness follows if we prove the estimate
\begin{equation}\label{thm2-goal2}
\|T_{\sigma}(a,g)\|_{L^p} \lesssim \|g\|_{L^{\infty}}
\end{equation}
for all $H^p$-atoms $a$.
Moreover, by virtue of the translation invariance,
\[
T_{\sigma}(a,g)(x+x_0)
=T_{\sigma_{x_0}}(a(\cdot+x_0), g(\cdot+x_0))(x),
\]
where $\sigma_{x_0}(x,\xi,\eta)=\sigma(x+x_0,\xi,\eta)$,
it is sufficient to treat $H^p$-atoms supported
in cubes centered at the origin.

Let $g \in L^{\infty}$ and 
let $a$ be an $H^p$-atom satisfying \eqref{atomQ} with a cube $Q$ 
centered at the origin and with $L>n/p-n$. 
We divide the $p$-th power of the $L^p$-norm
in the left hand side of \eqref{thm2-goal2} into
\begin{equation}\label{thm2-goal2-divide}
\|T_{\sigma}(a,g)\|_{L^p(Q^*)}^p
+\|T_{\sigma}(a,g)\|_{L^p((Q^*)^c)}^p.
\end{equation}
For the former term,
it follows from Theorem B with $(p,q,r)=(2,\infty,2)$ that
\[
\|T_{\sigma}(a,g)\|_{L^p(Q^*)}
\le |Q^*|^{1/p-1/2}\|T_{\sigma}(a,g)\|_{L^2}
\lesssim |Q|^{1/p-1/2}\|a\|_{L^2}\|g\|_{L^{\infty}}
\le \|g\|_{L^{\infty}},
\]
where we used the fact
\[
BS^{-(1-\rho)n/p}_{\rho,\rho} \subset 
BS^{-(1-\rho)n/2}_{\rho,\rho}.
\]
In the rest of this section,
we shall estimate the latter term in \eqref{thm2-goal2-divide}. 
The method will be similar to the one used 
in Section \ref{section3}.

Let $\Psi_j$, $j \ge 0$, be as in \eqref{dyadic-dec} with $d=2n$.
This time we do not need a delicate decomposition
such as \eqref{goal1-split-sum} for the proof of Theorem \ref{main-thm-2} 
and decompose $\sigma$ as
\[
\sigma(x,\xi,\eta)=\sum_{j \ge 0}\sigma_{j}(x,\xi,\eta)
\]
with
\[
\sigma_j(x,\xi,\eta)
=\sigma(x,\xi,\eta)\Psi_j(\xi,\eta).
\]
We also write the partial inverse Fourier transform of 
$\sigma_{j}(x,\xi, \eta)$ with respect to $(\xi, \eta)$ as
\[
K_{j}(x,y,z)
=\frac{1}{(2\pi)^{2n}}
\int_{(\R^n)^2}e^{i(y\cdot\xi+z\cdot\eta)}
\sigma_{j}(x,\xi,\eta)\, d\xi d\eta,
\quad
x,y,z \in \R^n,
\]
and then
\[
T_{\sigma_{j}}(a,g)(x)
=\int_{(\R^n)^2}
K_{j}(x,x-y,x-z)a(y)g(z)\, dydz.
\]
Notice that $\sigma_{j}$ satisfies the following: 
\begin{align}
& \supp \sigma_{j} (x, \cdot, \cdot) 
\subset 
\{ |\xi|\le 2^{j+1}, \; |\eta|\le 2^{j+1}\}, 
\label{sigmaj-1}
\\
& 
1+ |\xi|+ |\eta| \approx 2^{j} \quad \text{on}\quad 
\supp \sigma_{j} (x, \cdot, \cdot), 
\label{sigmaj-2}
\\
&
|\partial_{x}^{\alpha}\partial_{\xi}^{\beta}\partial_{\eta}^{\gamma}
\sigma_{j} (x, \xi, \eta) |
\lesssim 
2^{j m}\, 2^{j \rho (|\alpha|- |\beta| - |\gamma|)}. 
\label{sigmaj-3}
\end{align}

Let $x \not\in Q^*$.
By the same argument using the moment condition on $a$
as in \eqref{goal1-moment} with $c_Q=0$,
we have
\begin{equation}\label{goal2-moment}
\begin{split}
T_{\sigma_j}(a,g)(x)
&=L\sum_{|\alpha|=L}
\int_{\substack{y\in Q \\ z \in \R^n \\ 0<t<1}}
(1-t)^{L-1}\frac{(-y)^{\alpha}}{\alpha!}K_{j}^{(\alpha,0)}(x,x-ty,x-z)
\\
&\qquad \qquad \times  a(y)g(z)\, dydzdt, 
\end{split}
\end{equation}
where $K_{j}^{(\alpha,0)}(x,y,z)=\partial_{y}^{\alpha}K_{j}(x,y,z)$. 
We take integers $M$ and $M'$ satisfying
$M>n/p-n/2$ and $M'>n/2$.
Since $|x| \approx |x-ty|$
for $x \not\in Q^*$, $y \in Q$ and $0<t<1$,
Schwarz's inequality with respect to the $z$-variable gives 
\begin{align*}
&(1+2^{j\rho}|x|)^{M}
|T_{\sigma_{j}}(a,g)(x)|
\\
&\lesssim
\|g\|_{L^{\infty}}\ell(Q)^{L-n/p}
\sum_{|\alpha|=L}
\int_{\substack{y\in Q \\ z \in \R^n \\ 0<t<1}}
\\
&\qquad \times (1+2^{j\rho}|x-ty|)^{M}
|K_{j}^{(\alpha,0)}(x,x-ty,z)|
\, dydzdt
\\
&\lesssim 
\|g\|_{L^{\infty}}2^{-j\rho n/2}\ell(Q)^{L-n/p}
\sum_{|\alpha|=L}
\int_{\substack{y\in Q \\ 0<t<1}}
\\
&\qquad  \times
\left\|(1+2^{j\rho}|x-ty|)^{M}(1+2^{j\rho}|z|)^{M'}
K_{j}^{(\alpha,0)}(x,x-ty,z)\right\|_{L^2_z}
\, dydt.
\end{align*}
Thus, by writing
\begin{align*}
h_{j}^{(Q,L)}(x)
&=2^{-j\rho n/2}\ell(Q)^{L-n/p}
\sum_{|\alpha|=L}\sum_{|\beta| \le M}\sum_{|\gamma| \le M'}
\int_{\substack{y\in Q \\ 0<t<1}}
\\
&\qquad  \times
\left\|(2^{j\rho}(x-ty))^{\beta}(2^{j\rho}z)^{\gamma}
K_{j}^{(\alpha,0)}(x,x-ty,z)\right\|_{L^2_z}\, dydt,
\end{align*}
we have
\begin{equation}\label{goal2-estimate-moment}
|T_{\sigma_{j}}(a,g)(x)|
\lesssim
(1+2^{j\rho}|x|)^{-M}
h_{j}^{(Q,L)}(x)
\|g\|_{L^{\infty}}.
\end{equation}

We shall make a slight modification
on the argument \eqref{goal1-estimate-h_j^{(Q,L)}}
to estimate the $L^2$-norm of $h_{j}^{(Q,L)}$.
By Minkowski's inequality for integrals,
\begin{align*}
\|h_{j}^{(Q,L)}\|_{L^2}
&\le
2^{-j\rho n/2}\ell(Q)^{L-n/p}
\sum_{|\alpha|=L}\sum_{|\beta| \le M}\sum_{|\gamma| \le M'}
\int_{\substack{y\in Q \\ 0<t<1}}
\\
&\quad  \times
\left\|(2^{j\rho}(x-ty))^{\beta}(2^{j\rho}z)^{\gamma}
K_{j}^{(\alpha,0)}(x,x-ty,z)\right\|_{L^2_{x,z}}\, dydt.
\end{align*}
In the same way as in Section \ref{section3}, 
the function in the above $\|\cdot\|_{L^2_{x,z}}$
can be written as
\begin{align*}
&(2^{j\rho}(x-ty))^{\beta}(2^{j\rho}z)^{\gamma}
K_{j}^{(\alpha,0)}(x,x-ty,z)
\\
&=C_{\alpha,\beta,\gamma}
\int_{(\R^n)^2}
e^{i\{x\cdot\xi+z\cdot\eta\}}
(2^{j\rho}\partial_{\xi})^{\beta}(2^{j\rho}\partial_{\eta})^{\gamma}
\left[\xi^{\alpha}\sigma_{j}(x,\xi,\eta)\right]
\\
&\qquad \qquad \qquad \times
e^{-ity \cdot \xi}\psi_0(\xi/2^{j+1})\psi_0(\eta/2^{j+1})
\, d\xi d\eta,
\end{align*}
where $\psi_0$ is as in \eqref{dyadic-dec} with $d=n$
and we used \eqref{sigmaj-1}. 
From \eqref{sigmaj-2} and \eqref{sigmaj-3}, 
we see that
\[
\left|\partial_x^{\alpha'}\partial_{\xi}^{\beta'}\partial_{\eta}^{\gamma'}
(2^{j\rho}\partial_{\xi})^{\beta}(2^{j\rho}\partial_{\eta})^{\gamma}
\left[ \xi^{\alpha}\sigma_{j}(x,\xi,\eta)\right] 
\right|
\lesssim 2^{j(|\alpha|+m )}
(1+|\xi|+|\eta|)^{\rho(|\alpha'|-|\beta'|-|\gamma'|)}.
\]
Hence, the Calder\'on-Vaillancourt theorem on $\R^{2n}$ (\cite{CV}) 
and Plancherel's theorem give
\begin{align*}
&\left\|(2^{j\rho}(x-ty))^{\beta}(2^{j\rho}z)^{\gamma}
K_{j}^{(\alpha,0)}(x,x-ty,z)\right\|_{L^2_{x,z}}
\\
&\lesssim
2^{j(|\alpha|+m )}
\left\|e^{-ity \cdot \xi}\psi_0(\xi/2^{j+1})\psi_0(\eta/2^{j+1})
\right\|_{L^2_{\xi,\eta}}
\\
&\approx  2^{j(|\alpha|+m +n )},
\end{align*}
which implies
\begin{equation}\label{goal2-estimate-h_j^{(Q,L)}}
\begin{split}
\|h_j^{(Q,L)}\|_{L^2}
&\lesssim
2^{-j\rho n/2}\ell(Q)^{L-n/p}
\sum_{|\alpha|=L}
2^{j(|\alpha|+m+ n)}|Q|
\\
&\approx \left(2^j \ell(Q)\right)^{L-n/p+n}2^{j\rho n(1/p-1/2)}.
\end{split}
\end{equation}

If we do not use the moment condition on $a$ in \eqref{goal2-moment},
a similar argument yields
\begin{equation}\label{goal2-estimate-not-moment}
|T_{\sigma_{j}}(a,g)(x)|
\lesssim
(1+2^{j\rho}|x|)^{-M}
h_{j}^{(Q,0)}(x)
\|g\|_{L^{\infty}} 
\end{equation}
with
\begin{align*}
h_{j}^{(Q,0)}(x)
&=2^{-j\rho n/2}\ell(Q)^{-n/p}
\sum_{|\beta| \le M}\sum_{|\gamma| \le M'}
\int_{y\in Q}
\\
&\qquad \times
\left\|(2^{j\rho}(x-y))^{\beta}(2^{j\rho}z)^{\gamma}
K_{j}(x,x-y,z)\right\|_{L^2_z}\, dy
\end{align*}
and
\begin{equation}\label{goal2-estimate-h_j^{(Q,0)}}
\|h_{j}^{(Q,0)}\|_{L^2}
\lesssim
\left(2^j \ell(Q)\right)^{-n/p+n}2^{j\rho n(1/p-1/2)}.
\end{equation}

Combining \eqref{goal2-estimate-moment}
and \eqref{goal2-estimate-not-moment},
we have
\[
|T_{\sigma_j}(a,g)(x)|
\lesssim (1+2^{j\rho}|x|)^{-M}
\min\left\{h_j^{(Q,L)}(x), h_j^{(Q,0}(x)\right\}
\|g\|_{L^{\infty}}.
\]
Using \eqref{goal2-estimate-h_j^{(Q,L)}}, 
\eqref{goal2-estimate-h_j^{(Q,0)}} 
and H\"older's inequality with $1/p=1/q+1/2$, we have 
\begin{align*}
&\left\|(1+2^{j\rho}|\cdot|)^{-M}
\min\left\{h_j^{(Q,L)}, h_j^{(Q,0}\right\}\right\|_{L^p}
\\
&\lesssim \left\|(1+2^{j\rho}|\cdot|)^{-M}\right\|_{L^q}
\min\left\{\left\|h_j^{(Q,L)}\right\|_{L^2}, 
\left\|h_j^{(Q,0}\right\|_{L^2}\right\}
\\
&\lesssim
\min\left\{\left(2^j \ell(Q)\right)^{L-n/p+n},
\left(2^j \ell(Q)\right)^{-n/p+n}\right\}.
\end{align*}
Therefore,
\begin{align*}
&\|T_{\sigma}(a,g)\|_{L^p((Q^*)^c)}^p
\le \sum_{j=0}^{\infty} 
\|T_{\sigma _{j}}(a,g)\|_{L^p((Q^*)^c)}^p
\\
&\lesssim
\left(\sum_{2^j\ell(Q) \le 1}\left(2^j \ell(Q)\right)^{(L-n/p+n)p}
+\sum_{2^j\ell(Q)>1}\left(2^j \ell(Q)\right)^{(-n/p+n)p}\right)
\|g\|_{L^{\infty}}^{p} 
\lesssim
\|g\|_{L^{\infty}}^{p}, 
\end{align*}
which is the desired estimate
for the latter term in \eqref{thm2-goal2-divide}.


\end{document}